\pdfoutput=1
\documentclass[letterpaper, 10pt, conference]{ieeeconf} 

\IEEEoverridecommandlockouts                              
\usepackage{empheq,amssymb}
\usepackage{braket}
\usepackage{dsfont,mathrsfs}
\usepackage{enumerate}
\usepackage{graphicx}
\usepackage{amsmath}
\usepackage{stfloats}

\newtheorem{definition} {Definition}

\newtheorem{theorem}    {Theorem}
\newtheorem{lemma}      {Lemma}

\newcommand\EXP[1]{\mathop{\kern0pt \mathds E}{\Set{#1}}}
\newcommand\PR [1]{\mathop{\kern0pt \Pr}{\Set{#1}}}

\newcommand{\leftexp}[2]%
  {\mathop{}%
   \mathopen{\vphantom{#2}}^{#1}%
   \kern-\scriptspace%
   #2}

\makeatletter
\newcommand\SEQ{\@ifstar\SEQB\SEQA}
\makeatother

\newcommand\SEQA[2][T]{\{#2_t$, $t=1,\dots,#1\}}
\newcommand\SEQB[1]{\{#1_1$, $t=1,\dots\}}

\newcommand{\ind}{\mathds{1}}

\title {\LARGE \bf Decentralized Detection with Signaling}
\author{
Ashutosh Nayyar and Demosthenis Teneketzis 
\thanks{A.\ Nayyar and D.\ Teneketzis are with Department of Electrical Engineering and Computer Science,
        University of Michigan, Ann arbor. 
         {\tt\small anayyar@umich.edu, teneket@eecs.umich.edu}}
         }
\begin{document}
\maketitle
\begin{abstract}
    We consider a sequential problem in decentralized detection. Two observers can make repeated noisy observations of a binary hypothesis on the state of the environment. At any time, any of the two observers can stop and send a final message to the other observer or it may continue to take more measurements. After an observer has sent its final message, it stops operating.  The other observer is then faced 
with a different stopping problem. At each time instant, it can decide either to stop and declare a final decision
on the hypothesis or take another measurement. At each time, the system incurs an operating cost depending on the number of observers that are active at that time. A terminal cost that measures the accuracy of the final decision is incurred at the end. We show that, unlike in other sequential detection problems, stopping rules characterized by two thresholds on an observer's posterior belief no longer guarantee optimality in this problem. Thus the potential for signaling among observers alters the nature of optimal policies. We obtain a new parametric characterization of optimal policies for this problem.
\end{abstract}

\section{Introduction}
     Decentralized detection problems are motivated by applications in large scale decentralized systems such as sensor networks, power systems and surveillance networks. In such networks, sensors receive different information about the environment but share a common objective, for example to decide if a fault has occurred or not in a power system, or to detect the presence of a target in a surveillance area. Sensors may be allowed to communicate but they are constrained to exchange only a limited amount of information because of energy constraints, data storage and data processing constraints, communication constraints etc. 
     \par
    Decentralized detection problems may be static or sequential. In static problems, sensors make a fixed number of observations about a hypothesis on the state of the environment which is modeled as a random variable \(H\). Sensors may transmit a single message (a quantized version of their observations) to a fusion center which makes a final decision on \(H\). Such problems have been extensively studied since their initial formulation in \cite{Tenney_Detection} (See the surveys in \cite{Tsitsiklis_survey}, \cite{Varshney} and references therein).  In most such formulations, it has been shown that person-by-person optimal decision rules (as defined in \cite{Radner_team}) for a binary hypothesis detection problem are characterized by thresholds on the likelihood ratio (or equivalently on the posterior belief on the hypothesis). 
    \par
 
 In sequential problems, the number of observations taken by the sensors is
not fixed a priori. In the centralized sequential detection problems, as formulated in \cite{Wald}, a sensor can sequentially make costly observations  and, after each observation, can choose whether to stop and declare its final decision on $H$ or to take more observations. In the decentralized analogue of the sequential detection problem, two (or more) sensors locally decide when to stop taking more measurements and then make a final decision on $H$. Each sensor pays a penalty for delaying its final decision and a terminal cost that depends on the final decisions of all the sensors and the true value of $H$ is incurred at the end.
       A version of this problem (called the decentralized Wald problem) was formulated in
      \cite{Dec_Wald} and it was shown that at each time instant, optimal policies for the
      sensors are described by two thresholds. The computation of these thresholds requires solution of two coupled sets of dynamic programming equations. Similar results were obtained in a continuous time setting in \cite{LaVigna_86}. 
      \par
      A key feature of the decentralized Wald problem is that the individual sensors do not communicate their decisions to each other. That is, the $i^{th}$ sensor is not aware of decisions of other sensors. This implies that if policies of all other sensors are fixed, the $i^{th}$ sensor is faced with a classical sequential detection problem for which two-threshold policies are optimal. 
      In the problem we consider in this paper, each sensor observes the other sensor's decisions. Hence, in addition to its own measurements of $H$, the $i^{th}$ sensor can use the decisions made by other sensors (whether they have stopped or not and whether the final decision was $0$ or $1$) to make its decisions. The final decision of the sensor that stops in the end is taken as the final decision made by the group of sensors. Thus, sensors can convey information to each other through their decisions. The presence of signaling among sensors implies that, even if all other sensors have fixed their strategies, the problem for $i^{th}$ sensor is no longer a classical sequential detection problem. We show that, for this problem, the classical two-threshold policies no longer guarantee optimality. We obtain an alternative parametric characterization of the optimal policies of the sensors. A related sequential detection problem with one-way communication was presented in \cite{NayyarTeneketzis:2009}.
\par
 
         The rest of the paper is organized as follows. In Section~\ref{sec:PF}, we formulate our problem with two observers. We present the information states for the observers in Section~\ref{sec:IS}. A counterexample that shows that classical two-thresholds are not necessarily optimal is presented in Section~\ref{sec:counter}. We derive a parametric characterization of optimal policies in Section~\ref{sec:thresholds}. We conclude in Section~\ref{sec:con}.

\emph{Notation:} Throughout this paper, \(X_{1:t}\) refers to the sequence \(X_1, X_2,..,X_t\). Subscripts are used as time index and the superscripts are used as the index of the sensor. We use capital letters to denote random variable and the corresponding lower case letters for their realizations.

\section{Problem formulation} \label{sec:PF}
   Consider a binary hypothesis problem where the true hypothesis is modeled as a random variable \(H\) taking values 0 or 1 with known prior probabilities.
   \[P(H=0) = p_0; \hspace{10pt} P(H=1) = 1-p_0\]
 Consider two observers: observer 1 (O1) and observer 2 (O2). We assume that each observer can make noisy observations of the true hypothesis. Conditioned on the hypothesis \(H\), the following statements are assumed to be true: \\
 1. The observation of the \(i^{th}\) observer at time \(t\), \((Y^i_t)\) (taking values in the set \(\mathcal{Y}^i\)), either has a discrete distribution \((P^i_t(.|H))\) or admits a probability density function \((f^i_t(.|H))\). \\
2. Observations of the \(i^{th}\) observer at different time instants are conditionally independent given \(H\).\\
3. The observation sequences at the two observers are conditionally independent given \(H\). \\
\begin{figure}[ht]
\begin{center}

\includegraphics[height=5cm,width=7cm]{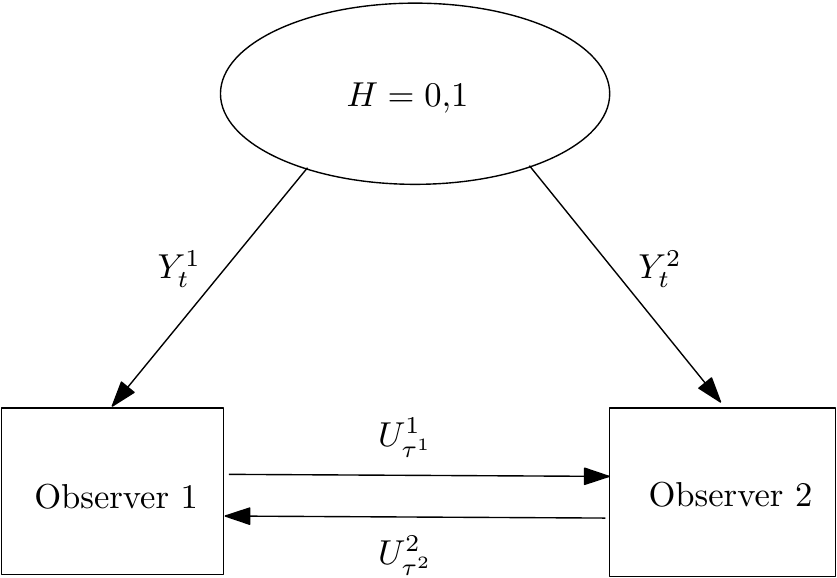}
\caption{Decentralized Detection}

\end{center} 
\end{figure}
Observer $i$ ($i=1,2$) observes the measurement process \(Y^i_t\)\((t=1,2,...)\). If no observer has stopped before time \(t\), then at time $t$, any observer can decide either to stop and send a binary message 0 or 1 to the other observer or to postpone its decision and get another measurement. After an observer has sent its final message, it stops operating.  The other observer (the one which has not yet stopped) is then faced
with a different stopping problem. At each time instant, it can decide either to stop and declare a final decision
on the hypothesis or take another measurement. 

  We denote by $U^i_t$ the decision of the $i^{th}$ observer at time $t$. $U^i_t$ belongs to the set $\{0,1,b\}$, where $b$ denotes a blank, that is, no message or no final decision. At time $t$, the $i^{th}$ observer makes its decision based on its observations till time $t$ and the messages (blanks or otherwise) exchanged between the two observers till time $t-1$. We have,
 \begin{equation}
 U^i_t = \gamma^i_t(Y^i_{1:t},U^1_{1:t-1},U^2_{1:t-1})
 \end{equation}
 where the collection of functions $\Gamma^i:=(\gamma^i_t$, $t=1,2,\ldots)$ constitute the policy of the $i^{th}$ observer. We define the following stopping times:
\[\tau^1 := min\{t:U^1_t \neq b\}\]
\[\tau^2 := min\{t:U^2_t \neq b\}\]
After the time $\tau^i$, all future observations and decisions of the $i^{th}$ observer are assumed to be in the empty set ($\emptyset$), that is, the observer stops taking measurements or making decisions. We assume that the final decision on the hypothesis must be made no later than a finite horizon $T$, hence we have that $\tau^1 \leq T$ and $\tau^2 \leq T$. We define, $ ~~\tau^{min} := min\{\tau^1,\tau^2\}$, $\tau^{max} := max\{\tau^1,\tau^2\} $ and
\[L:= \left\{ \begin{array}{l} 1 \mbox{~~~if~~~} \tau^2 < \tau^1 \\
                               2 \mbox{~~~if~~~} \tau^2 \geq \tau^1 \end{array} \right. \]
The system incurs a total cost given by:
 \begin{align}
  &\mathcal{J}(\Gamma^1, \Gamma^2) \nonumber \\ 
  &:= \mathds{E}^{\Gamma^1,\Gamma^2}\{K(\tau^{min}-1) + k(\tau^{max}-\tau^{min}) \nonumber\\
  &~~~~~~~~~~~~~+ J(U^{L}_{\tau^{max}},H)\} \label{eq:objective}
 \end{align}    
 where $K>k>0$ are constants and $J(\cdot,\cdot)$ is a non-negative distortion function with $J(0,0)=J(1,1)=0$. (The superscript $\Gamma^1,\Gamma^2$ over the expectation denotes that the expectation is with respect to a measure that depends on the choice of the policies $\Gamma^1,\Gamma^2$.) The first term in the objective represents the operating costs when both observers are active, the second term represents the operating cost of only one observer. The operating costs incorporate the cost of taking a new measurement, the energy cost of staying on for another time step and/or a penalty for delaying the decision. The last term in the objective represents the accuracy of the final decision. Note that only the decision of the observer that stops later is considered as the final decision. In case of simultaneous decisions, only observer~2's decision is considered as the final decision. We can now formulate the optimization problem as follows:
 \medskip 
 \\
\emph{Problem P:} Given the statistics of the binary hypothesis and the observation processes, the cost parameters $K,k$, the distortion function $J(\cdot,\cdot)$ and a time horizon $T$, the objective is to select policies $\Gamma^1,\Gamma^2$ that minimize the total expected cost in (\ref{eq:objective}).

\section{Information States} \label{sec:IS}
\begin{figure*}[!b]
\line(1,0){500}
\begin{align}
 V_t(\pi,0) = min \{&\mathds{E}[\ind_{\{\tau^1=t\}}J(0,H)+ \ind_{\{\tau^1>t\}}(k(\tau^1-t)+ J(U^1_{\tau^1},H))|\Pi^2_{t}=\pi,\ind_{\{\tau^1<t\}=0},U^2_t=0], \nonumber \\
 &\mathds{E}[\ind_{\{\tau^1=t\}}J(1,H) + \ind_{\{\tau^1>t\}}(k(\tau^1-t)+ J(U^1_{\tau^1},H))|\Pi^2_{t}=\pi,\ind_{\{\tau^1<t\}=0},U^2_t=1], \nonumber \\
 &\mathds{E}[\ind_{\{\tau^1=t\}}(k+V_{t+1}(\Pi^2_{t+1},1)) \nonumber \\
 &~~+\ind_{\{\tau^1>t\}}(K+V_{t+1}(\Pi^2_{t+1},0))|\Pi^2_{t}=\pi,\ind_{\{\tau^1<t\}=0},U^2_t=b] \label{eq:defv0}
 \}
\end{align}   
\end{figure*}  

In this section, we identify information states for the two observers. We start by fixing the policy of observer~1 to an arbitrary choice and finding a sufficient statistic for observer~2. The nature of this sufficient statistic does not depend on the arbitrary choice of observer~1's policy. This sufficient statistic is the information state for observer 2.
\par
Consider a fixed policy \(\Gamma^1 = (\gamma^1_1,\gamma^1_2,...,\gamma^1_{T})\) for O1. At any time $t$, we define the following:
\medskip
\begin{definition}   
  Given a fixed policy $\Gamma^1$ of observer~1 and functions $\gamma^2_{1:t-1}$, we define observer~2's belief on the hypothesis given all its information at time $t$. \[ \Pi^2_t := P(H=0|Y^2_{1:t},U^1_{1:t-1},U^2_{1:t-1}=b_{1:t-1}), \] 
 where $b_{1:t-1}$ denotes a sequence of blank messages from time $1$ to $t-1$. For \(t=0\), we have \(\Pi^2_0 = p_0\).\\
    

 \end{definition}
\medskip
We will show that the pair $(\Pi^2_t, \ind_{\{\tau^1<t\}})$ is the information state for observer~2. We first describe the evolution of $(\Pi^2_t, \ind_{\{\tau^1<t\}})$ in time in the following lemma.
\medskip
\begin{lemma}\label{lemma:update}
\begin{enumerate}[(i)]
\item $\ind_{\{\tau^1<t+1\}} = \ind_{\{\tau^1<t\}} + \ind_{\{U^1_t \neq b\}}$\\
\item  With observer~1's policy fixed to $\Gamma^1$, $\Pi^2_t$ evolves as follows:\\  $ \Pi^2_{t+1} = \left \{ \begin{array}{ll}
               f_{t+1}(\Pi^2_t,Y^2_{t+1}) & \mbox{if $\ind_{\{\tau^1<t\}} =1$} \\
               g_{t+1}(\Pi^2_t,Y^2_{t+1},U^1_{t}) & \mbox{if $\ind_{\{\tau^1<t\}} =0$}
               \end{array}
               \right.$, \\
~~~~\\
where $f_{t+1}$ and $g_{t+1}$ are deterministic functions.

\end{enumerate}
\end{lemma}   
\begin{proof} See Appendix \ref{sec:lemma_1}. \end{proof}
\medskip
The optimal policy for observer~2 (for the given choice of $\Gamma^1$) can be obtained by means of a dynamic program. We now define the value functions of the dynamic program.
\medskip
\begin{definition}\label{def:value_functions}
\begin{enumerate}[(i)]
\item For $\pi \in [0,1]$ and $a \in \{0,1\}$, we define 
\begin{align}
   V_T(\pi,a) := min \{
   &\mathds{E}[J(0,H)|\Pi^2_{T}=\pi], \nonumber \\
  	                                  &\mathds{E}[J(1,H)|\Pi^2_{T}=\pi]\}\nonumber
\end{align}
\item For $\pi \in [0,1]$ and $t=T-1,T-2,\ldots,1$, we define
\begin{align}
V_t(\pi,1)& \nonumber \\= min &\{
               \mathds{E}[J(0,H)|\Pi^2_{t}=\pi],
  	                                 \mathds{E}[J(1,H)|\Pi^2_{t}=\pi], \nonumber \\
 	                                  &k+ \mathds{E}[V_{t+1}(\Pi^2_{t+1},1)|\Pi^2_{t}=\pi,\ind_{\{\tau^1<t\}=1}] \} \nonumber
\end{align} 
\item  For $\pi \in [0,1]$ and $t=T-1,T-2,\ldots,1$, we define $V_t(\pi,0)$ in equation (\ref{eq:defv0}) at the bottom of the page.
\end{enumerate}
\end{definition}       
\medskip
\begin{theorem}\label{thm:DP}
  With a fixed policy \(\Gamma^1\) of observer~1, there is an optimal policy for observer~2 of the form:
  \begin{equation*} \label{eq:Qp1}
  U^2_t = \gamma^2_{t}(\Pi^{2}_{t},\ind_{\{\tau^1<t\}})
  \end{equation*}
  for \(t=1,2,...,T\). Moreover, this optimal policy can be obtained by the dynamic program described by the value functions in Definition~\ref{def:value_functions}. Thus, at time $t$ and for a given $\pi$ and $a$, the optimal decision is $0$ (or $1$/$b$) if the first (or second/third) term is the minimum in the definition of $V_t(\pi,a)$.
\end{theorem}
\begin{proof} See Appendix \ref{sec:info_states_proof}.
\end{proof}
\medskip
The above arguments can be repeated by interchanging the roles of observer~1 and observer~2 to conclude that for a fixed policy $\Gamma^2$ of observer~2, an optimal policy for observer~1 is of the form:
\begin{equation*} \label{eq:Qp2}
  U^1_t = \gamma^1_{t}(\Pi^{1}_{t},\ind_{\{\tau^2<t\}})
\end{equation*}
 where   $\Pi^1_t := P^{\Gamma^2}(H=0|Y^1_{1:t},U^2_{1:t-1},U^1_{1:t-1}=b_{1:t-1})$. Note, however, that the actual dynamic program for observer~1 will differ from that of Theorem~\ref{thm:DP} because of the asymmetry in the objective function of equation~(\ref{eq:objective}) when both observers make simultaneous decisions to stop. While the value functions $V_T$ and $V_t(\pi,1)$ in the dynamic program for observer~1 will the be same as in Definition 2, the value function $V_t(\pi,0)$ is given in equation (\ref{eq:defv0forO1}) at the bottom of the next page.     
\section{A Counterexample} \label{sec:counter}
\begin{figure*}[!b]
\line(1,0){500}
\begin{align}
 V_t(\pi,0) = min \{&\mathds{E}[k(\tau^2-t)+ J(U^2_{\tau^2},H)|\Pi^1_{t}=\pi,\ind_{\{\tau^2<t\}=0},U^1_t=0], \nonumber \\
 &\mathds{E}[k(\tau^2-t)+ J(U^2_{\tau^2},H)|\Pi^1_{t}=\pi,\ind_{\{\tau^2<t\}=0},U^1_t=1], \nonumber \\
 &\mathds{E}[\ind_{\{\tau^2=t\}}(k+V_{t+1}(\Pi^1_{t+1},1)) \nonumber \\
 &~~+\ind_{\{\tau^2>t\}}(K+V_{t+1}(\Pi^1_{t+1},0))|\Pi^1_{t}=\pi,\ind_{\{\tau^2<t\}=0},U^1_t=b] \label{eq:defv0forO1}
 \}
\end{align}   
\end{figure*}  
In the sequential detection problem with a single observer \cite{Wald} , it is well known that an optimal policy is a function of the observer's posterior belief \(\Pi_t\) and is described by two thresholds at each time \(t\). That is the decision at time \(t\), \(Z_t\), is given as:
      \[ Z_t = \left \{ \begin{array}{ll}
               1 & \mbox{if $\Pi_t \leq \alpha_t$} \\
               b & \mbox{if $\alpha_t<\Pi_t < \beta_t$} \\
               0 & \mbox{if $\Pi_t \geq \beta_t$}
               \end{array}
               \right. \]
     where \(b\) denotes a decision to continue taking measurement and \(\alpha_t \leq \beta_t\) are real numbers in \([0,1]\).
 A similar two-threshold structure of optimal policies was also established for the decentralized Wald problem in \cite{Dec_Wald}. We will show by means of a counterexample that such a structure is not necessarily optimal in our problem.
 \par

 Consider the following instance of Problem P. We have equal prior on \(H\), that is \(P(H=0)=P(H=1)=1/2\) and a time horizon of \(T=3\). Assume $k=1$ and $1<K<2$. The observation space of observer~1 is \(\mathcal{Y}^1 = \{0,1\}\) and the observations at time \(t\) obey the following conditional probabilities:
 \[ \begin{array}{lcc}
Observation & 0 & 1 \\
P(\cdot|H=0)    & q_t & (1-q_t)  \\
P(\cdot|H=1)    & (1-q_t) & q_t  
\end{array}\]
where $q_1 =q_2 =1/2$ and $q_3=1$. Thus, the first two observations of observer~1 reveal no information about $H$ while the third observation reveals $H$ noiselessly.
The observation space of observer~2 is \(\mathcal{Y}^1 = \{0,1,2\}\) and the observations at time \(t\) obey the following conditional probabilities:
\[ \begin{array}{lccc}
Observation & 0 & 1 & 2 \\
P(\cdot|H=0)    & r_t & (1-r_t) & 0 \\
P(\cdot|H=1)    & 0 & (1-r_t) & r_t 
\end{array}\]
where \(r_2=r_3=0\) and $0<r_1<1$. Thus, the second and third observations of observer~2 reveal no information about $H$. Note that under this statistical model of observations, there exists a choice of policies such that the system makes perfect final decision on the hypothesis and incurs only operational costs (if observer~2 stops at $t=1$ and observer~1 waits till time $t=3$, then it can make a perfect decision on $H$ and the system incurs an operational cost of $2k =2$). We assume that the cost of  a mistake in the final decision \((U^L_{\tau_{max}} \neq H)\) is sufficiently high so that any choice of policies that makes a mistake in the final decision with non-zero probability will have a performance worse than $2$. Thus, any choice of policies that makes a mistake in the final decision with non-zero probability cannot be optimal.
\par

Under the above instance of our problem, if observer~2 is restricted to use a two-threshold rule at time $t=1$, then the lowest achievable value of the objective is given as:
\begin{align} \label{eq:threshold_cost}
\min[\{r_1 + (1-r_1)(K+1)\}, \{2-r_1/2\}] \end{align}
where the first term corresponds to the case when $\gamma^2_1$ is given as:
\begin{equation} \label{eq:ex_1}
 U^2_1 = \left \{ \begin{array}{ll}
               1 & \mbox{if $\Pi^2_1 = 0$} \\
               b & \mbox{if $0<\Pi^2_1<1$} \\
               0 & \mbox{if $\Pi^2_1 = 1$}
               \end{array}
               \right. \end{equation}
  and the second term corresponds to $\gamma^2_1$ being
\begin{equation} \label{eq:ex_2}  
 U^2_1 = \left \{ \begin{array}{ll}
               1 & \mbox{if $\Pi^2_1 < 1$} \\
               0 & \mbox{if $\Pi^2_1 = 1$}
               \end{array}
               \right. \end{equation}
 Other choices of thresholds for observer~2 at time $t=1$ do not give a lower value than the expression in (\ref{eq:threshold_cost}).
 \par
 Consider now the following choice of $\gamma^{2,*}_1$:
 \begin{equation}\label{eq:non_th}
  U^2_1 = \left \{ \begin{array}{ll}
               1 & \mbox{if $\Pi^2_1 = 0$} \\
               0 & \mbox{if $0<\Pi^2_1<1$} \\
               b & \mbox{if $\Pi^2_1 = 1$}
               \end{array}
               \right. \end{equation}
 The lowest achievable expected cost under the above choice of $\gamma^2_1$ is $\mathcal{J}^{*} = 2(1-r_1) + r_1(K+1)/2$. It is easy to check that for $1<K<2$ and $r_1<2/3$,
 \[ \mathcal{J}^{*} < \min[\{r_1 + (1-r_1)(K+1)\}, \{2-r_1/2\}]\]
 Thus, $\gamma^{2,*}_1$ outperforms the two-threshold rules.  
 \medskip
 
 \emph{Discussion:} In the above example, observer~1 can always make the correct decision at time $t=3$. However, this incurs additional operational costs. A good policy should try to enable the observers to make the correct decision before time $t=3$, whenever possible. If observer~2 gets the observations \(0\) or \(2\) at \(t=1\), then it is certain about the true hypothesis. Using the first threshold rule given in (\ref{eq:ex_1}), observer 2 is able to convey to observer~1 that it is certain about the true hypothesis and what this hypothesis is, thus preventing observer~1 from waiting till time $t$ $=$ $3$ to make the final decision.  However, in the case when observer~2 gets measurement \(1\), it sends a blank and postpones its decision to stop. This incurs additional operating costs for observer~2 without providing any new information that may prevent observer~1 from delaying its decision to time $t=3$. By making $r_1$ small, the contribution of this term in the overall cost is increased. The second threshold (equation (\ref{eq:ex_2})) rule attempts to keep the operational costs of observer~2 small but does not always send enough information to observer~1 to enable it to make a decision before $t=3$ even when observer~1 knows the true value of $H$. The non-threshold rule (given in (\ref{eq:non_th})), however, keeps the operational cost of observer~2 small (when $r_1$ is small) while at the same time ensuring that whenever observer~2 is certain about true value of $H$, the final decision is not postponed to time $t=3$. 
\section{Parametric Characterization of Optimal Policies} \label{sec:thresholds}

An important advantage of the threshold rules in the case of the centralized or the decentralized Wald problem is that it modifies the problem of finding the globally optimal policies from a sequential functional optimization problem to a sequential parametric optimization problem. Even though we have established that a classical threshold rule does not hold for our problem, it is still possible to get a finite parametric characterization of optimal policies. Such a parametric characterization provides significant computational advantage in finding optimal policies, for example by reducing the search space for an optimal policy.         \par
 In Theorem~\ref{thm:DP}, we have established  that for an arbitrarily fixed choice of observer~1's policy, the optimal policy for observer~2 can be determined by a dynamic program using the value functions \(V_t(\pi,a), t=T,...,2,1\). 
 We have the following lemma.
 \begin{lemma} \label{lemma:lemma_2}
  With a fixed (but arbitrary) choice of \(\Gamma^1\), the value function at \(T\) can be expressed as:
  \begin{equation}
  V_{T}(\pi,a) := min\{ l^{0}(\pi), l^{1}(\pi)\}
  \end{equation}
  where \(l^{0}\) and \(l^{1}\) are affine functions of \(\pi\).
  Also, the value functions at time \(t\) can be expressed as:
  \begin{equation}\label{eq:thresholds_lemma1}
  V_{t}(\pi,1) := min\{ l^{0}(\pi), l^{1}(\pi), G_t(\pi)\}
  \end{equation}
  where  \(G_t\) is a concave function of \(\pi\), and
  \begin{equation}\label{eq:thresholds_lemma2}
  V_{t}(\pi,0) := min\{ L^{0}_t(\pi), L^{1}_t(\pi), H_t(\pi)\}
  \end{equation}
  where \(L^{0}_{t}\) and \(L^{1}_{t}\) are affine functions of \(\pi\) and \(H_t\) is a concave function of \(\pi\) (the actual form of these functions depends on the choice of \(\Gamma^1\)).
   \end{lemma}
 \begin{proof} See Appendix \ref{sec:lemma_2}. \end{proof}
 \begin{theorem}
  For any fixed policy \(\Gamma^1\) of observer~1,  an optimal policy for observer~2 can be characterized as follows:
   \[ U^2_{T} = \left \{ \begin{array}{ll}
               1 & \mbox{if $\Pi^2_{T} \leq \alpha_{T}$} \\
               0 & \mbox{if $\Pi^2_{T} > \alpha_{T}$ }
               \end{array}
               \right. \]
               where \(0 \leq \alpha_{T} \leq 1\). For \(t=1,2,..,T^1-1\), if $\ind_{\{\tau^1<t\}} = 1$,
               \[ U^2_{t} = \left \{ \begin{array}{ll}
               1 & \mbox{if $\Pi^2_{t} \leq \alpha_{t}(1)$} \\
               b & \mbox{if $\alpha_t(1) < \Pi^2_{t} \leq \beta_{t}(1)$} \\
               0 & \mbox{if $\Pi^2_{t} > \beta_{t}(1)$ }
               \end{array}
               \right. \]
               where \(0 \leq \alpha_t(1) \leq \beta_t(1) \leq 1\), 
               and if $\ind_{\{\tau^1<t\}} = 0$,
         \[ U^2_t = \left \{ \begin{array}{ll}
               b & \mbox{if $\Pi^2_t < \alpha_t(0)$}  \\
               1 & \mbox{if $\alpha_t(0) \leq \Pi^2_t \leq \beta_t(0)$} \\
               b & \mbox{if $\beta_t(0)<\Pi^2_t < \delta_t(0)$} \\
               0 & \mbox{if $\delta_t(0) \leq \Pi^2_t \leq \theta_t(0)$} \\
               b & \mbox{if $\Pi^2_t > \theta_t(0)$}
               \end{array}
               \right. \]  
               where \(0 \leq \alpha_t(0) \leq \beta_t(0) \leq \delta_t(0) \leq \theta_t(0) \leq 1\). 
\end{theorem}
\begin{proof}
From lemma~\ref{lemma:lemma_2}, we know that the value functions can be written as minimum of affine and concave functions. Since taking minimum of two straight lines and a concave function can partition the interval \( [0,1]\) into at most five intervals, this gives a four threshold characterization of optimal policy where the thresholds signify the boundaries of these intervals. At time $T$ or when $\ind_{\{\tau^1<t\}} = 1$, observer~2's decision of $0$ or $1$ is the final decision ($U^{L}_{\tau_{max}}$) on $H$. In these cases, if observer~2 is certain about $H$ (that is its belief is $0$ or $1$), then it should clearly choose the correct value of $H$. This fact reduces the number of thresholds for time $T$ and when $\ind_{\{\tau^1<t\}} = 1$. 
\end{proof}
\medskip

\emph{Discussion:} It is instructive to compare our problem with the decentralized Wald problem studied in \cite{Dec_Wald}. Both problems involve two observers that make repeated observations of $H$ and decide when to stop. Unlike the problem in this paper, in the decentralized Wald problem the sensors do not have access to each other's decisions, that is, an observer's policy is restricted to be of the form:
\[ U^i_t = \gamma^i_t(Y^i_{1:t},U^i_{1:t-1}) \]
The optimality of classical two-threshold rules for the decentralized Wald problem was established in \cite{Dec_Wald}. In this paper, we allowed each observer to observe other's decisions and showed that the two threshold rules are no longer optimal.
\par
Both the decentralized Wald problem and the problem formulated in this paper are \emph{team problems}. That is, they involve more than one decision maker with a common objective. However, in the decentralized Wald problem, a decision-maker's decisions do not influence the information available to other decision-makers. This is the essential criterion for \emph{static} team problems \cite{Ho1980}. The problem formulated in this paper is a \emph{dynamic} team problem since a decision-maker's past decisions are a part of the information available to other decision-makers. It is the dynamic aspect of this team problem that allows for signaling between decision-makers and changes the nature of optimal policies.   
\section{Conclusions}\label{sec:con}
   We considered a sequential problem in decentralized detection problem with signaling. Two observers make separate costly measurements of a binary hypothesis and decide when to stop. The observers can observe each other's decisions (whether the other observe has stopped or not and whether the final decision was $0$ or $1$). The final decision of the observer that stops in the end is taken as the final decision made by the group. Thus, observers can convey information to each other through their decisions. We identified information states for the two observers and showed that classical two  threshold rules no longer guarantee optimality. However, a finite parametric characterization of optimal policies is still possible. 

\appendices

 \section{Proof of Lemma 1} \label{sec:lemma_1}
 \begin{proof}
 Part (i) follows from definition of $\tau^1$.
  
 In Part (ii), if $\tau^1 < t$, then by definition, we have 
 \begin{align} \Pi^2_t &:= P(H=0|Y^2_{1:t},U^1_{1:t-1},U^2_{1:t-1}=b_{1:t-1}) \notag \\
 &= P(H=0|Y^2_{1:t},U^1_{1:\tau^1}) \label{eq:mtnsapp1.4}
  \end{align}
 where we removed redundant terms from the conditioning (terms which are constants or functions of other terms).
 Similarly, 
 \begin{align} \Pi^2_{t+1} =  P(H=0|Y^2_{1:t+1},U^1_{1:\tau^1}), \notag  
 \end{align}
 which, on using Bayes' rule gives,
 \begin{align}
   \Pi^2_{t+1}&= \frac{P(Y^2_{t+1}|H=0)\Pi^2_t}{P(Y^2_{t+1}|H=0)\Pi^2_t + P(Y^2_{t+1}|H=1)(1-\Pi^2_t)} \notag \\
   &=: f_{t+1}(\Pi^2_{t+1},Y^2_{t+1})\label{eq:Ap1}
 \end{align}
 
 If $\tau^1 \geq t$, then $U^1_{1:t-1}=b_{1:t-1}$ (that is all decisions of observer 1 are blanks till time $t-1$) and
 \begin{align} \Pi^2_{t} &:= P(H=0|Y^2_{1:t},U^1_{1:t-1}=b_{1:t-1},U^2_{1:t-1}=b_{1:t-1}) \notag 
 \end{align}
 Also,
 \begin{align} &\Pi^2_{t+1} := P(H=0|Y^2_{1:t+1},U^1_{1:t-1}=b_{1:t-1},U^1_t,U^2_{1:t}=b_{1:t}) \notag \\
 &=\frac{P(Y^2_{t+1},U^1_t,H=0|Y^2_{1:t},U^1_{1:t-1}=b_{1:t-1},U^2_{1:t}=b_{1:t})}{\displaystyle\sum_{h\in \{0,1\}}P(Y^2_{t+1},U^1_t,h|Y^2_{1:t},U^1_{1:t-1}=b_{1:t-1},U^2_{1:t}=b_{1:t})} \label{eq:mtnsapp1.1}
 \end{align}
 The numerator in (\ref{eq:mtnsapp1.1}) can be written as:
 \begin{align}
 &P(Y^2_{t+1}|H=0)\cdot \notag \\&\{P(U^1_t|H=0,Y^2_{1:t},U^1_{1:t-1}=b_{1:t-1},U^2_{1:t}=b_{1:t})\}\Pi^2_t \label{eq:mtnsapp1.3}
 \end{align}
  
 We now focus on the second term in (\ref{eq:mtnsapp1.3}).\\
 \emph{Claim:} Consider a realization $u^1_t$, $y^2_{1:t}$. Then, 
 \begin{align}
 &P(u^1_t|H=0,y^2_{1:t},U^1_{1:t-1}=b_{1:t-1},U^2_{1:t}=b_{1:t}) \notag \\
 &= P(u^1_t|H=0,U^1_{1:t-1}=b_{1:t-1},U^2_{1:t}=b_{1:t}) \label{eq:mtnsapp1.7}
 \end{align}
 Moreover, under the given choice of $\Gamma^1$, the probability on the right hand side of (\ref{eq:mtnsapp1.7}) is a function only of $u^1_t$.
 \par
 \emph{Proof of claim:} Using Bayes' rule,
 \begin{align}
 &P(u^1_t|H=0,y^2_{1:t},U^1_{1:t-1}=b_{1:t-1},U^2_{1:t}=b_{1:t}) \notag \\
 &= \frac{P(u^1_t,H=0,y^2_{1:t},U^1_{1:t-1}=b_{1:t-1},U^2_{1:t}=b_{1:t})}{\displaystyle\sum_{u'}P(U^1_t=u',H=0,y^2_{1:t},U^1_{1:t-1}=b_{1:t-1},U^2_{1:t}=b_{1:t})} \label{eq:mtnsapp1.10}
 \end{align}
 Consider the joint probability in the numerator in (\ref{eq:mtnsapp1.10})
 \begin{align}
 P(u^1_t,H=0,y^2_{1:t},b^1_{1:t-1},b^2_{1:t}) \notag 
 \end{align}
 where we use $b^1_{1:t-1},b^2_{1:t-1}$ as shorthand notations for  $U^1_{1:t-1}=b_{1:t-1}$ and $U^2_{1:t-1}=b_{1:t-1}$ respectively. This probability can be further written as:  
 \begin{align}
 &=\displaystyle\sum_{y^1_{1:t}}P(u^1_t,H=0,y^2_{1:t},b^1_{1:t-1},b^2_{1:t},y^1_{1:t}) \notag\\ 
 &=\displaystyle\sum_{y^1_{1:t}}[ P(u^1_t|y^1_{1:t},b^1_{1:t-1},b^2_{1:t-1})\notag\\&\cdot P(U^2_t=b|y^2_{1:t},b^1_{1:t-1},b^2_{1:t-1}) \notag\\
 &\cdot P(y^1_t|H=0)P(y^2_t|H=0)\notag\\ &\cdot \displaystyle\prod_{k=1}^{t-1}\{P(U^1_k=b|y^1_{1:k},b^1_{1:k-1},b^2_{1:k-1})\notag\\&\cdot P(U^2_k=b|y^2_{1:k},b^1_{1:k-1},b^2_{1:k-1}) \notag\\&\cdot P(y^1_k|H=0)P(y^2_k|H=0)\} ] \cdot p_0 \label{eq:mtnsapp1.5}
 \end{align}
 Rearranging the summation in (\ref{eq:mtnsapp1.5}), we get
 \begin{align}
 &P(U^2_t=b|y^2_{1:t},b^1_{1:t-1},b^2_{1:t-1})P(y^2_t|H=0) \cdot p_0 \notag\\
 &\cdot\displaystyle\prod_{k=1}^{t-1}\{P(U^2_k=b|y^2_{1:k},b^1_{1:k-1},b^2_{1:k-1}) P(y^2_k|H=0)\} \notag \\
 &\cdot\displaystyle\sum_{y^1_{1:t}}[ P(u^1_t|y^1_{1:t},b^1_{1:t-1},b^2_{1:t-1})P(y^1_t|H=0)\notag\\ &\cdot \displaystyle\prod_{k=1}^{t-1}\{P(U^1_k=b|y^1_{1:k},b^1_{1:k-1},b^2_{1:k-1}) P(y^1_k|H=0)\} ]  \label{eq:mtnsapp1.6}
 \end{align}
 Expressions similar to (\ref{eq:mtnsapp1.6}) hold for each term in the denominator of (\ref{eq:mtnsapp1.10}) and the terms outside the summation over $y^1_{1:t}$ cancel out in the numerator and the denominator. We note that the summation over $y^1_{1:t}$ in (\ref{eq:mtnsapp1.6}) does not depend on $y^2_{1:t}$. Hence, the conditional probability in the left hand side of (\ref{eq:mtnsapp1.10}) does not depend on $y^2_{1:t}$. This establishes equation (\ref{eq:mtnsapp1.7}). We also note that under the fixed policy $\Gamma^1$ of observer~1, the summation over $y^1_{1:t}$ in (\ref{eq:mtnsapp1.6}) is a function only of $u^1_t$. Thus, the probability on the right hand side of (\ref{eq:mtnsapp1.7}) is a function only of $u^1_t$. This concludes the proof of the claim.
\par
 
Using the result of the claim in (\ref{eq:mtnsapp1.3}) and using similar arguments for the denominator in (\ref{eq:mtnsapp1.1}), we get  
 \begin{align} &\Pi^2_{t+1} \notag\\
 &=\frac{P(Y^2_{t+1}|H=0)P(U^1_t|H=0,b^1_{1:t-1},b^2_{1:t})\Pi^2_t}{\begin{array}{l}P(Y^2_{t+1}|H=0)P(U^1_t|H=0,b^1_{1:t-1},b^2_{1:t})\Pi^2_t \\ +P(Y^2_{t+1}|H=1)P(U^1_t|H=1,b^1_{1:t-1},b^2_{1:t})(1-\Pi^2_t) \end{array} } \notag \\
 &=: g_{t+1}(\Pi^2_{t+1},Y^2_{t+1},U^1_t) \label{eq:g_def}
 \end{align}
 \end{proof}
\section{Proof Outline of Theorem 1}\label{sec:info_states_proof}
We provide an outline of the proof of Theorem~1. The general idea is to show that at each time $t$, the value functions of Definition~2 represent the optimal future costs. Therefore, a policy that for each realization of $\Pi^2_t, \ind_{\{\tau^1<t\}}$ selects the minimizing term in the corresponding value function achieves the optimal cost. Thus, an optimal policy can be found that depends only on $\Pi^2_t, \ind_{\{\tau^1<t\}}$. We start from time $T$.

\par
 If the observer~2 is active at the terminal time $T$, it can only make one of two decisions: $0$ or $1$. The expected future cost of choosing $u \in \{0,1\}$ for observer~2 is
\begin{align}
&E[J(u,H)|Y^2_{1:T},U^1_{1:T-1},U^2_{1:T-1}=b_{1:T-1}] \notag \\
&=  J(u,0)\Pi^2_T + J(u,1)(1-\Pi^2_T)\notag\\
&= E[J(u,H)|\Pi^2_{T}] \label{eq:app2.2}
\end{align}
Thus, the value function at time $T$ is the minimum of the expected future costs incurred by choosing $0$ or $1$. Hence, it represents the optimal expected future cost for observer~2 at time $T$. Proceeding backwards, we assume that the value functions at time $t+1,t+2,\ldots,T$ represent optimal future costs at the respective times and consider two cases at each time $t<T$. 
\medskip

\emph{Case A: $\tau^1<t$} If observer~1 has already stopped before $t$, then observer~2's stopping problem is the same as the centralized Wald problem and the value function $V_t(\pi,1)$ is same as the value function in the dynamic program for the Wald problem.

\emph{Case B: $\tau^1 \geq t$} We now consider the case when observer~1 has not stopped before time $t$. If observer~2 decides to stop and chooses $u \in \{0,1\}$ at time $t$, then the expected future cost will be
\begin{align}
 \mathds{E}\Big[&\ind_{\{\tau^1=t\}}J(u,H) \notag \\&+ \ind_{\{\tau^1>t\}}(k(\tau^1-t)+ J(U^1_{\tau^1},H))\bigg|\begin{array}{l}Y^2_{1:t},b^1_{1:t-1}\\,b^2_{1:t-1},U^2_t=u \end{array}\Big] \label{eq:app2.1}
\end{align}

\emph{Claim:} The expectation in (\ref{eq:app2.1}) is same as:
\begin{align}
&\mathds{E}[\ind_{\{\tau^1=t\}}J(u,H) + \notag \\&\ind_{\{\tau^1>t\}}(k(\tau^1-t)+ J(U^1_{\tau^1},H))|\Pi^2_t,\ind_{\{\tau^1<t\}}=0,U^2_t=u]
\end{align}
\emph{Proof of claim:}
For each realization $y^2_{1:t}$ of observer~2's observations, the expectation in (\ref{eq:app2.1}) depends on the conditional distribution of the following random variables: $H,\tau^1,U^1_{\tau^1}$ given the realization of the random variables $y^2_{1:t},b^1_{1:t-1},b^2_{1:t-1},U^2_t=u$. Note that under the fixed policy $\Gamma^1$ of observer~1, $\tau^1,U^1_{\tau^1}$ are functions of observer~1's observation sequence $Y^1_{1:T}$ and the terms $b^1_{1:t-1},b^2_{1:t-1},U^2_t=u$ fixed in the conditioning.  Hence the conditional belief
\[ P(H,\tau^1,U^1_{\tau^1}|y^2_{1:t},b^1_{1:t-1},b^2_{1:t-1},U^2_t=u) \]
is a deterministic transformation of the belief
\[ P(H,Y^1_{1:T}|y^2_{1:t},b^1_{1:t-1},b^2_{1:t-1},U^2_t=u).\]
We will show that the above probability is same as 
\[ P(H,Y^1_{1:T}|\pi^2_t,\ind_{\{\tau^1<t\}} =0,U^2_t=u)\]
and hence the conditional expectation in (\ref{eq:app2.1}) is same as 
\begin{align*} \mathds{E}[&\ind_{\{\tau^1=t\}}J(u,H)+ \notag \\&\ind_{\{\tau^1>t\}}(k(\tau^1-t)+ J(U^1_{\tau^1},H))|\pi^2_{t},\ind_{\{\tau^1<t\}=0},U^2_t=u]\end{align*}
which corresponds to the first two terms in the minimization in $V_t(\pi^2_t,0)$ in equation (\ref{eq:defv0}).

Consider $P(H=0,y^1_{1:T}|y^2_{1:t},b^1_{1:t-1},b^2_{1:t-1},U^2_t=u)$
\begin{align}
&=P(y^1_{1:T}|H=0,y^2_{1:t},b^1_{1:t-1},b^2_{1:t-1},U^2_t=u)\pi^2_t \label{eq:app2.5} 
\end{align}
Similarly, 
\begin{align}
&P(H=0,y^1_{1:T}|\pi^2_t,\ind_{\{\tau^1<t\}}=0,U^2_t=u) \notag \\
&= P(y^1_{1:T}|H=0,\pi^2_t,\ind_{\{\tau^1<t\}}=0,U^2_t=u)\pi^2_t \label{eq:app2.6}
\end{align}
We now compare the first terms in (\ref{eq:app2.5}) and (\ref{eq:app2.6}). Consider the first term in (\ref{eq:app2.5}), which can be written as
\begin{align}
&= \frac{P(H=0,y^1_{1:T},y^2_{1:t},b^1_{1:t-1},b^2_{1:t-1},U^2_t=u)}{\displaystyle\sum_{\tilde{y}^1_{1:T}} P(H=0,\tilde{y}^1_{1:T},y^2_{1:t},b^1_{1:t-1},b^2_{1:t-1},U^2_t=u)} \label{eq:app2.3}
\end{align}
The numerator can be written as:
\begin{align}
 &P(y^1_{t+1:T}|H=0)P(U^2_t=u|y^2_{1:t},b^1_{1:t-1},b^2_{1:t-1}) \notag\\
 &\cdot P(y^1_t|H=0)P(y^2_t|H=0)\notag\\ &\cdot \displaystyle\prod_{k=1}^{t-1}\{P(U^1_k=b|y^1_{1:k},b^1_{1:k-1},b^2_{1:k-1})\notag\\&\cdot P(U^2_k=b|y^2_{1:k},b^1_{1:k-1},b^2_{1:k-1}) \notag\\&\cdot P(y^1_k|H=h)P(y^2_k|H=h)\} \cdot p_0
\end{align}
Similar expressions hold for the denominator in (\ref{eq:app2.3}) and the terms that depend on $y^2_{1:t}$ will cancel in the numerator and the denominator. Therefore,
\begin{align}
&P(y^1_{1:T}|H=0,y^2_{1:t},b^1_{1:t-1},b^2_{1:t-1},U^2_t=u) \notag\\
&= P(y^1_{1:T}|H=0,b^1_{1:t-1},b^2_{1:t-1},U^2_t=u)\label{eq:app2.7} 
\end{align}
Now consider the first term in (\ref{eq:app2.6}) which can be written as:
\begin{align}
&P(y^1_{1:T}|H=0,\pi^2_t,\ind_{\{\tau^1<t\}}=0,U^2_t=u) \notag\\
&=\displaystyle \sum_{y^2_{1:t}}[P(y^1_{1:T}|y^2_{1:t},H=0,\pi^2_t,\ind_{\{\tau^1<t\}}=0,U^2_t=u) \notag\\
&\cdot P(y^2_{1:t}|H=0,\pi^2_t,\ind_{\{\tau^1<t\}}=0,U^2_t=u)] \notag \\
&= \displaystyle \sum_{y^2_{1:t}}[P(y^1_{1:T}|y^2_{1:t},H=0,b^1_{1:t-1},b^2_{1:t-1},U^2_t=u) \notag\\
&\cdot P(y^2_{1:t}|H=0,\pi^2_t,\ind_{\{\tau^1<t\}}=0,U^2_t=u)] \label{eq:neweq1}
\end{align}
The first term inside the summation in (\ref{eq:neweq1}) is same as LHS of (\ref{eq:app2.7}). Using (\ref{eq:app2.7}) in (\ref{eq:neweq1}) gives
\begin{align}
& \displaystyle \sum_{y^2_{1:t}}[P(y^1_{1:T}|H=0,b^1_{1:t-1},b^2_{1:t-1},U^2_t=u) \notag\\
&\cdot P(y^2_{1:t}|H=0,\pi^2_t,\ind_{\{\tau^1<t\}}=0,U^2_t=u)] \notag \\
&= P(y^1_{1:T}|H=0,b^1_{1:t-1},b^2_{1:t-1},U^2_t=u) \label{eq:app2.8}
\end{align} 
which is same as RHS of (\ref{eq:app2.7}). Thus the probabilities in RHS of (\ref{eq:app2.5}) and (\ref{eq:app2.6}) are equal. Similar conclusions hold for $H=1$ in (\ref{eq:app2.5}) and (\ref{eq:app2.6}) . This implies  the equality of expectations and completes the proof of the claim.
\par
As a consequence of the claim, the first two terms in the minimization in the definition of $V_t(\pi,0)$ (equation \ref{eq:defv0}) correspond to the expected future cost of choosing $0$ or $1$ at time $t$. On the other hand, if observer~2 decides to continue at time $t$, then by the fact that value functions at $t+1$ represent the expected future costs at $t+1$, we can write the expected future cost as:
\begin{align}
&\mathds{E}\Big[\ind_{\{\tau^1=t\}}(k+V_{t+1}(\Pi_{t+1},1)) \nonumber \\
 &~~+\ind_{\{\tau^1>t\}}(K+V_{t+1}(\Pi_{t+1},0))\bigg|\begin{array}{l}Y^2_{1:t},b^1_{1:t-1}\\,b^2_{1:t-1},U^2_t=b \end{array}\Big] \label{eq:app2.1.B}
\end{align}

Using lemma 1 and the fact that under fixed policy $\Gamma^1$ of observer~1, $U^1_{t}$ is a function of observer~1's observation sequence $Y^1_{1:t}$ and the terms $b^1_{1:t-1},b^2_{1:t-1}$ fixed in the conditioning , one can conclude that for each realization of $y^2_{1:t}$ this expectation in (\ref{eq:app2.1.B}) is a function of the following conditional probability:\[P(Y^1_{1:t},Y^2_{t+1}|y^2_{1:t},b^1_{1:t-1},b^2_{1:t-1},U^2_t=b)\]
Using arguments similar to those in the claim above, it can be shown that the above conditional probability is same as:
\[ P(Y^1_{1:t},Y^2_{t+1}|\pi^2_t,\ind_{\{\tau^1<t\}} =0,U^2_t=b)\]
This shows that the third term in the minimization in the definition of $V_t(\pi,0)$ (equation \ref{eq:defv0}) is the expected future cost of making a decision to continue at time $t$. Thus,
$V_t(\Pi^2_t,0)$ is the minimum of the future costs incurred by choosing $0,1$ or $b$. Hence, it represents the optimal future cost at time $t$, if observer~1 has not already stopped before time $t$.

\section{Proof Outline of Lemma 2} \label{sec:lemma_2}
 We define the following functions
 \begin{align*}
 l^0(\pi) :&= J(0,0)\pi+J(0,1)(1-\pi) \notag\\
 &=\mathds{E}[J(0,H)|\Pi^2_t =\pi] \\
  l^1(\pi) :&= J(1,0)\pi+J(1,1)(1-\pi) \notag\\
 &=\mathds{E}[J(1,H)|\Pi^2_t =\pi]
 \end{align*} 
 For the value function at time $T$, the result of the lemma follows from the definitions of $l^0(\pi), l^1(\pi)$ and $V_T(\pi,a)$. Since, for each $a \in \{0,1\}$, $V_T(\pi,a)$ is the minimum of two affine functions of $\pi$, it implies that, for each $a \in \{0,1\}$, $V_T(\pi,a)$ is a concave function of $\pi$. Now, assume that $V_{t+1}(\pi,a)$ is concave in $\pi$ for each $a \in \{0,1\}$. The concavity of the value functions at time $t+1$ implies that they can be written as infimum of affine functions of $\pi$. In particular, we have
 \begin{align} \label{eq:concave_1}
 V_{t+1}(\pi,1) = \inf_{i}\{ a_i\pi+b_i\}
 \end{align} 
 and 
\begin{align} \label{eq:concave_2}
 V_{t+1}(\pi,0) = \inf_{i} \{c_i\pi+d_i\}
 \end{align}
 Now consider $V_t(\pi,1)$. The first two terms in the definition of $V_t(\pi,1)$ are affine in $\pi$ (see Definition 2). We need to show that the third term -
 \begin{equation} \label{eq:app3.1}
 k+ \mathds{E}[V_{t+1}(\Pi^2_{t+1},1)|\Pi^2_{t}=\pi,\ind_{\{\tau^1<t\}=1}] \end{equation}
 -is a concave function of $\pi$. From equation (\ref{eq:Ap1}) in the proof of lemma 1, we know that $\Pi_{t+1}$ can be written as:
 \begin{align}
 \Pi^2_{t+1}&= \frac{P(Y^2_{t+1}|H=0)\Pi^2_t}{P(Y^2_{t+1}|H=0)\Pi^2_t + P(Y^2_{t+1}|H=1)(1-\Pi^2_t)} \notag \\
            &= \frac{P(Y^2_{t+1}|H=0)\Pi^2_t}{P(Y^2_{t+1}|\Pi^2_t)}  \label{eq:app3.2}
 \end{align}
 Substituting (\ref{eq:app3.2}) in (\ref{eq:app3.1}) and evaluating the expectation gives:
 \begin{align}
  k + \displaystyle\sum_{y \in \mathcal{Y}}&P(Y^2_{t+1}=y|\Pi^2_{t}=\pi,\ind_{\{\tau^1<t\}}=1) \notag\\
                                           &V_{t+1}\left(\frac{P(Y^2_{t+1}=y|H=0)\pi}{P(Y^2_{t+1}=y|\pi)},1\right) 
 \end{align}
 Now using the characterization of $V_{t+1}(\pi,1)$ from (\ref{eq:concave_1}), we get
 \begin{align}
  k + \displaystyle\sum_{y \in \mathcal{Y}}&P(Y^2_{t+1}=y|\Pi^2_{t}=\pi) \notag\\
                                           &\left[ \inf_i\{ a_i\left(\frac{P(Y^2_{t+1}=y|H=0)\pi}{P(Y^2_{t+1}=y|\pi)}\right) +b_i \} \right] \notag\\
  =k + \displaystyle\sum_{y \in \mathcal{Y}}  & \inf_i \{ a_i(P(Y^2_{t+1}=y|H=0)\pi) \notag \\&+b_iP(Y^2_{t+1}=y|H=0)\pi\notag\\&+ b_iP(Y^2_{t+1}=y|H=1)(1-\pi) \}                                \end{align}
 Each term in the summation over $y \in \mathcal{Y}$ is infimum of affine functions of $\pi$, hence each term in the summation is a concave function of $\pi$. Thus,  the third term of $V_t(\pi,1)$ is a concave function of $\pi$.
  
  \par
  Next consider $V_t(\pi,0)$ defined in (\ref{eq:defv0}). The conditional expectation for the first (or second) term in the minimization in RHS of (\ref{eq:defv0}) is an affine function of the conditional probability of the random variables $H, Y^1_{1:T}$. Using arguments from Appendix~\ref{sec:info_states_proof}, this conditional probability can be written as:
  \begin{align}
  &P(y^1_{1:T},H=0|\pi^2_t,b^1_{1:t-1},b^2_{1:t-1},U^2_t=0) \notag \\
  &= P(y^1_{1:T}|H=0,b^1_{1:t-1},b^2_{1:t-1},U^2_t=0)\pi^2_t
 \end{align}
  and
  \begin{align}
  &P(y^1_{1:T},H=1|\pi^2_t,b^1_{1:t-1},b^2_{1:t-1},U^2_t=0) \notag \\
  &= P(y^1_{1:T}|H=1,b^1_{1:t-1},b^2_{1:t-1},U^2_t=0)(1-\pi^2_t)
 \end{align}
 Thus, the conditional probability $P(Y^1_{1:T},H|\pi^2_t,b^1_{1:t-1},b^2_{1:t-1},U^2_t=0)$ is an affine function of $\pi^2_t$. This establishes the affine nature of the first two terms of RHS of (\ref{eq:defv0}) . The third term in $V_t(\pi,0)$ can be written as:
 \begin{align}
 &\mathds{E}\Big[\ind_{\{U^1_t \neq b\}}(k+V_{t+1}(g_{t+1}(\pi,Y^2_{t+1},U^1_t),1)) \nonumber \\
 &+\ind_{\{U^1_t=b\}}(K+V_{t+1}(g_{t+1}(\pi,Y^2_{t+1},b),0))\Big|\begin{array}{l}\Pi^2_{t}=\pi,\\\ind_{\{\tau^1<t\}=0},\\U^2_t=b\end{array}\Big] \label{eq:app3.3}
 \end{align}
 Consider the first term in the summation in (\ref{eq:app3.3}). Evaluating the expectation, we get
 \begin{align}
 &\displaystyle \sum_{u^1_t,y^2_{t+1}} [\ind_{\{u^1_t \neq b\}}(k+V_{t+1}(g_{t+1}(\pi,y^2_{t+1},u^1_t),1)) \notag \\ &\cdot P(u^1_t,y^2_{t+1}|\Pi^2_{t}=\pi,\ind_{\{\tau^1<t\}=0},U^2_t=b)]\label{eq:neweq2}
 \end{align}
  Using the characterization of $g_{t+1}$ from (\ref{eq:mtnsapp1.1}) and (\ref{eq:g_def}), the characterization of $V_{t+1}(\pi,1)$ from (\ref{eq:concave_2}) and arguments from Appendix II, (\ref{eq:neweq2}) can be shown to be equal to
  \begin{align}
 &\displaystyle \sum_{u^1_t=0,1}\displaystyle \sum_{y^2_{t+1}}[k \notag\\&+ \inf_i \{ a_iP(y^2_{t+1}|H=0)P(u^1_t|H=0,b^1_{1:t-1},b^2_{1:t})\pi \notag \\
 &+b_iP(y^2_{t+1}|H=0)P(u^1_t|H=0,b^1_{1:t-1},b^2_{1:t})\pi \notag\\&+ b_iP(y^2_{t+1}|H=1)P(u^1_t|H=1,b^1_{1:t-1},b^2_{1:t})(1-\pi)\}]
 \end{align}
 which is concave in $\pi$ (since it is infimum of affine functions in $\pi$). Similar arguments can be made for the second term in (\ref{eq:app3.3}) to conclude the concavity of third term in $V_t(\pi,0)$.
\section*{Acknowledgments}
  This research was supported in part by NSF Grant CCR-0325571 and NASA Grant NNX09AE91G.
\bibliographystyle{IEEEtran}
\bibliography{myref}
\end{document}